\documentclass[12pt]{article}
 \usepackage{latexsym}
 \usepackage{amssymb}
 \usepackage{graphicx}

 \newtheorem{theorem}{Theorem}

\newtheorem{Theorem}[theorem]{Theorem}

\newtheorem{Proposition}[theorem]{Proposition}
\newtheorem{Assumption[theorem]}{Assumption}
\newtheorem{Lemma}[theorem]{Lemma}
\newtheorem{Corollary}[theorem]{Corollary}

\newcommand{\T}{{\cal T}}

\newcommand{\A}{{\cal A}}
\newcommand{\B}{{\cal B}}
\newcommand{\C}{{\cal C}}
\newcommand{\D}{{\cal D}}
\newcommand{\J}{{\cal J}}
\newcommand{\LL}{{\cal L}}
\newcommand{\Q}{{\cal Q}}

\newcommand{\1}{{\bf 1}}
\newcommand{\ii}{{\bf e}^{(i)}}
\newcommand{\jj}{{\bf e}^{(j)}}

\newcommand{\qed}{\nobreak \ifvmode \relax \else
      \ifdim\lastskip<1.5em \hskip-\lastskip
      \hskip1.5em plus0em minus0.5em \fi \nobreak
      \vrule height0.75em width0.5em depth0.25em\fi}

\def \ep{\hbox{ }\hfill$\Box$}

\begin{document}
\title{H$^+$-Eigenvalues of Laplacian and Signless Laplacian Tensors}

\author{
Liqun Qi \thanks{Email: maqilq@polyu.edu.hk. Department of Applied
Mathematics, The Hong Kong Polytechnic University, Hung Hom,
Kowloon, Hong Kong. This author's work was supported by the Hong
Kong Research Grant Council (Grant No. PolyU 501909, 502510, 502111
and 501212).}}

\date{\today} \maketitle

%---------------------------------------------------------------------------------Abstract
\begin{abstract}
\noindent We propose a simple and natural definition for the
Laplacian and the signless Laplacian tensors of a uniform
hypergraph.  We study their H$^+$-eigenvalues, i.e., H-eigenvalues
with nonnegative H-eigenvectors, and H$^{++}$-eigenvalues, i.e.,
H-eigenvalues with positive H-eigenvectors.  We show that each of
the Laplacian tensor, the signless Laplacian tensor and the
adjacency tensor has at most one H$^{++}$-eigenvalue, but has
several other H$^+$-eigenvalues.  We identify their largest and
smallest H$^+$-eigenvalues, and establish some maximum and minimum
properties of these H$^+$-eigenvalues. We then define analytic
connectivity of a uniform hypergraph and discuss its application in
edge connectivity. \vspace{3mm}

\noindent {\bf Key words:}\hspace{2mm} Laplacian tensor, signless
Laplacian tensor, uniform hypergraph, H$^+$-eigenvalue \vspace{3mm}

\noindent {\bf AMS subject classifications (2010):}\hspace{2mm} 05C65; 15A18
  \vspace{3mm}

\end{abstract}

%------------------------------------------------------------------------------Section 1

\section{Introduction}
\hspace{4mm} Recently, several papers appeared on spectral
hypergraph theory via tensors \cite{BP, CD, HQ, LQY, Li, PZ, XC,
XC1, XC2}.  These works are all on uniform hypergraphs \cite{Be}. In
2008, Lim \cite{Li} proposed to study spectral hypergraph theory via
eigenvalues of tensors.  In 2009, Bul\`{o} and Pelillo \cite{BP}
gave new bounds on the clique number of a graph based on analysis of
the largest eigenvalue of the adjacency tensor of a uniform
hypergraph. In 2012, Hu and Qi \cite{HQ} proposed a definition for
the Laplacian tensor of an even uniform hypergraph, and analyzed its
connection with edge and vertex connectivity. In the same year,
Cooper and Dutle \cite{CD} analyzed the eigenvalues of the adjacency
tensor (hypermatrix) of a uniform hypergraph, and proved a number of
natural analogs of basic results in spectral graph theory.  Li, Qi
and Yu \cite{LQY} proposed another definition for the Laplacian
tensor of an even uniform hypergraph, established a variational
formula for its second smallest Z-eigenvalue, and used it to provide
lower bounds for the bipartition width of the hypergraph.   In
\cite{XC, XC2}, Xie and Chang proposed a definition for the signless
Laplacian tensor of an even uniform hypergraph, studied its largest
and smallest H-eigenvalues and Z-eigenvalues, and its applications
in the edge cut and the edge connectivity of the hypergraph.  They
also studied the largest and the smallest Z-eigenvalues of the
adjacency tensor of a uniform hypergraph in \cite{XC1}. In
\cite{PZ}, Pearson and Zhang studied the H-eigenvalues and the
Z-eigenvalues of the adjacency tensor of a uniform hypergraph.

Precisely speaking, the tensors mentioned above may be called
hypermatrices.   In physics and mechanics, tensors are physical
quantities, while hypermatrices are multi-dimensional arrays.   In
geometry, a tensor to a hypermatrix is like a linear transformation
to a matrix - the former objects are defined without choosing bases
\cite{Qi1}. However, for the most papers in tensor decomposition,
spectral theory of tensors and spectral hypergraph theory, as the
most papers cited in this paper, the word ``tensors'' are used for
those multi-dimensional arrays.   Following this habit, we use the
word ``tensors'' in this paper.

A uniform hypergraph is also called a $k$-graph \cite{Be, BH}. Let
$G = (V, E)$ be a $k$-graph, where $V = \{ 1, 2, \ldots, n \}$ is
the vertex set, $E = \{e_1, e_2, \ldots, e_m\}$ is the edge set,
$e_p \subset V$ and $|e_p| = k$ for $p = 1, \ldots, m$, and $k \ge
2$. If $k=2$, then $G$ is an ordinary graph. We assume that $e_p
\not = e_q$ if $p \not = q$.  Two vertices are called {\sl adjacent}
if they are in the same edge.   Two vertices $i$ and $j$ are called
{\sl connected} if either $i$ and $j$ are adjacent, or there are
vertices $i_1, \ldots, i_s$ such that $i$ and $i_1$, $i_k$ and $j$,
$i_r$ and $i_{r+1}$ for $r = 1, \ldots, s-1$, are adjacent
respectively.   A $k$-graph $G$ is called {\sl connected} if any
pair of its vertices are connected.   The {\bf adjacency tensor} $\A
= \A(G)$ of $G$, is a $k$th order $n$-dimensional symmetric tensor,
with $\A = \left(a_{i_1i_2\cdots i_k}\right)$, where
$a_{i_1i_2\cdots i_k} = {1 \over (k-1)!}$ if $(i_1, i_2, \ldots,
i_k) \in E$, and $0$ otherwise. Thus, $a_{i_1i_2\cdots i_k} = 0$ if
two of its indices are the same. For $i \in V$, its degree $d(i)$ is
defined as $d(i) = \left| \{ e_p : i \in e_p \in E \} \right|$.  We
assume that every vertex has at least one edge. Thus, $d(i) > 0$ for
all $i$.  The {\bf degree tensor} $\D = \D(G)$ of $G$, is a $k$th
order $n$-dimensional diagonal tensor, with its $i$th diagonal entry
as $d(i)$.  We denote the maximum degree, the minimum degree and the
average degree of $G$ by $\Delta$, $\delta$ and $\bar d$
respectively. If $\bar d = \Delta = d$, then $G$ is a regular graph,
called a $d$-regular $k$-graph.

The definition of the adjacency tensor is natural.  It was studied
in \cite{BP, CD, XC1}.  On the other hand, the definitions of
Laplacian and signless Laplacian tensors in \cite{HQ, LQY, XC, XC2}
are based upon some forms of sums of $k$-th powers.   They are not
simple and natural, and only work when $k$ is even.

In this paper, we propose a simple and natural definition for the
Laplacian and the signless Laplacian tensors of a $k$-graph $G$.
Recall that when $k=2$, the Laplacian matrix and the signless
Laplacian matrix of $G$ are defined as $\LL = \D - \A$ and $\Q = \D
+ \A$ \cite{BH}.  Many results of spectral graph theory are based
upon this definition.   Thus, for $k \ge 3$, we propose to define
the {\bf Laplacian tensor} and the {\bf signless Laplacian tensor}
of $G$ simply by $\LL = \D - \A$ and $\Q = \D + \A$.   This
definition is simple and natural, and is closely related with the
adjacency tensor $\A$. Furthermore, the signless Laplacian tensor
$\Q$ is a symmetric nonnegative tensor, while the Laplacian tensor
$\LL$ is the limit of symmetric $M$-tensors in the sense of
\cite{ZQZ}. $M$-tensors are closely related with nonnegative tensors
\cite{ZQZ}. Thus, we may use the recently developed theory and
algorithms on eigenvalues of nonnegative tensors \cite{CPZ, CPZ1,
FGH, HHQ, LZI, NQZ, Qi2, YY, ZQX} to study $\LL$ and $\Q$.

We discover that $\LL$ and $\Q$ have very nice spectral properties.
They are not irreducible in the sense of \cite{CPZ}.   But they are
weakly irreducible in the sense of \cite{FGH} if $G$ is connected.
When $k \ge 3$, each of them has at least $n+1$ H-eigenvalues with
nonnegative H-eigenvectors. We call such H-eigenvalues {\bf
H$^+$-eigenvalues}. Furthermore, each of them has at most one
H$^+$-eigenvalue with a positive eigenvector. We call such an
H$^+$-eigenvalue an {\bf H$^{++}$-eigenvalue}.

The remainder of this paper is distributed as follows.  In the next
section, we review the definition and properties of eigenvalues and
H-eigenvalues of tensors, and introduce H$^+$-eigenvalues and
H$^{++}$-eigenvalues.   We study H$^+$-eigenvalues of $\A$, $\LL$
and $\Q$ in Section 3.   We show that each of $\A$, $\LL$ and $\Q$
has at most one H$^{++}$-eigenvalue, but has several other
H$^+$-eigenvalues.  In Sections 4, we study the smallest
H-eigenvalue of $\LL$, and its link with connectedness of $G$. We
identify the largest H$^+$-eigenvalue of $\LL$, and establish a
maximum property of this H$^+$-eigenvalue in Section 5.   We
establish some maximum properties of the largest H-eigenvalues of
$\Q$ and $\A$, and discuss methods for computing them in Section 6.
In Section 7, we identify the smallest H$^+$-eigenvalue of $\Q$,
establish a minimum property of this H$^+$-eigenvalue, and discuss
its applications in edge connectivity and maximum cut. In Section 8,
we define {\bf analytic connectivity} of $G$ as a minimum quantity
related with $\LL$, and discuss its application in edge
connectivity.  Some final remarks are made in Section 9.

Denote by $\1$ the all $1$ $n$-dimensional vector, $\1_j = 1$ for
$j=1, \ldots, n$.   Denote by $\ii$ the $i$th unit vector in
$\Re^n$, i.e., $\ii_j = 1$ if $i=j$ and $\ii_j = 0$ if $i \not =j$,
for $i, j = 1, \ldots, n$.  For a vector $x$ in $\Re^n$, we define
its support as supp$(x) = \{ i \in V : x_i \not = 0 \}$. Denote the
set of all nonnegative vectors in $\Re^n$ by $\Re^n_+$ and the set
of all positive vectors in $\Re^n$ by $\Re^n_{++}$. For a $k$th
order $n$-dimensional tensor $\C = \left(c_{i_1 \cdots i_k}\right)$,
$|\C |$ is a $k$th order $n$-dimensional tensor $|\C| =
\left(|c_{i_1 \cdots i_k}|\right)$.   If both $\C = \left(c_{i_1
\cdots i_k}\right)$ and $\B = \left(b_{i_1 \cdots i_k}\right)$ are
real $k$th order $n$-dimensional tensors, and $b_{i_1 \cdots i_k}
\le c_{i_1 \cdots i_k}$ for $i_1, \ldots, j_k = 1, \ldots, n$, then
we write $\B \le \C$.  We use $\J$ to denote the $k$th order
$n$-dimensional tensor with all of its entries being $1$.

\section{H$^+$-Eigenvalues and H$^{++}$-Eigenvalues}

In this section, we will review the definition and properties of
eigenvalues and H-eigenvalues of tensors in \cite{Qi}, introduce
H$^+$-eigenvalues and H$^{++}$-eigenvalues, and review the
Perron-Frobenius Theorem for nonnegative tensors in \cite{CPZ, FGH,
YY}.   We also discuss the reducibility and weak irreducibility of
$\LL$ and $\Q$ in this section.

Consider a real $k$th order $n$-dimensional tensor $\T =
\left(t_{i_1 \cdots i_k}\right)$.   Let $x \in C^n$.  Then
$$\T x^k = \sum_{i_1, \ldots, i_k = 1}^n t_{i_1 \cdots i_k}x_{i_1}
\cdots x_{i_k},$$ and $\T x^{k-1}$ is a vector in $C^n$, with its
$i$th component defined by
$$\left(\T x^{k-1}\right)_i = \sum_{i_2, \ldots, i_k = 1}^n t_{i i_1 \cdots i_k}x_{i_2}
\cdots x_{i_k}.$$ Let $r$ be a positive integer.  Then $x^{[r]}$ is
a vector in $C^n$, with its $i$th component defined by $x_i^r$. We
say that $\T$ is symmetric if its entries $t_{i_1 \cdots i_k}$ are
invariant under any permutation of its indices.

Suppose that $x \in C^n$, $x \not = 0$, $\lambda \in C$,  $x$ and
$\lambda$ satisfy
\begin{equation} \label{eig}
\T x^{k-1} = \lambda x^{[k-1]}.
\end{equation}
Then we call $\lambda$ an {\bf eigenvalue} of $\T$, and $x$ its
corresponding {\bf eigenvector}.   From (\ref{eig}), we may see that
if $\lambda$ is an eigenvalue of $\T$ and $x$ is its corresponding
eigenvector, then
\begin{equation} \label{eig1}
\lambda = {(\T x^{k-1})_j \over x_j^{k-1}},
\end{equation}
for some $j$ with $x_j \not = 0$.   In particular, if $x$ is real,
then $\lambda$ is also real. In this case, we say that $\lambda$ is
an {\bf H-eigenvalue} of $\T$ and $x$ is its corresponding {\bf
H-eigenvector}.    If $x \in \Re^n_+$, then we say that $\lambda$ is
an {\bf H$^+$-eigenvalue} of $\T$.  If $x \in \Re^n_{++}$, then we
say that $\lambda$ is an {\bf H$^{++}$-eigenvalue} of $\T$.   If
$\lambda$ is an {\bf H$^+$-eigenvalue} but not an {\bf
H$^{++}$-eigenvalue} of $\T$, then we say that $\lambda$ is a {\bf
strict H$^+$-eigenvalue} of $\T$.

We say that $\T$ is positive definite (semi-definite) if $\T x^k >
0$ ($\T x^k \ge 0$) for all $x \in \Re^n, x \not = 0$.   Clearly,
$\T$ is positive definite only if $k$ is even, and when $k$ is odd,
$\T$ is positive semi-definite only if $\T$ is the zero tensor.

Note that (\ref{eig}) is a homogeneous system of $x$, with $n$
variables and $n$ equations.  We may regard that these variables
take values in the complex field.   According to algebraic geometry
\cite{CLO}, the {\sl resultant} of (\ref{eig}) is a polynomial in
the coefficients of (\ref{eig}), hence a polynomial in $\lambda$,
which vanishes if and only if (\ref{eig}) has a nonzero solution
$x$. Denote this polynomial by $\phi_\T(\lambda)$, and call it the
{\bf characteristic polynomial} of $\T$.

The main properties of eigenvalues and H-eigenvalues of a real $k$th
order $n$-dimensional symmetric tensor in \cite{Qi} are summarized
in the following theorem.

\begin{theorem} {\bf (Eigenvalues of Real Symmetric Tensors)} {\bf (Qi 2005)} \label{t1}

The followings hold for the eigenvalues of a real $k$th order
$n$-dimensional symmetric tensor $\T$:

(a). A number $\lambda \in C$ is an eigenvalue of $\T$ if and only
if it is a root of the characteristic polynomial $\phi_\T$.  Hence,
we regard the multiplicity of an eigenvalue $\lambda$ of $\T$ as its
multiplicity as a root of $\phi_\T$.

(b). The number of eigenvalues of $\T$, counting their
multiplicities, is $n(k-1)^{n-1}$. Their product is equal to
det$(\T)$, the resultant of $\T x^{k-1} = 0$.

(c). The sum of all the eigenvalues of $T$ is
$$(k-1)^{n-1}{\rm tr}(\T),$$
where tr$(\T)$ denotes the sum of the diagonal entries of $\T$.

(d). If $k$ is even, then $\T$ always has H-eigenvalues.   $\T$ is
positive definite (positive semi-definite) if and only if all of its
H-eigenvalues are positive (nonnegative).

(e). The eigenvalues of $\T$ lie in the following $n$ disks:

$$ | \lambda - t_{ii \cdots i} | \le \sum \left\{ | t_{ii_2
\cdots i_k} | : i_2, \ldots, i_k = 1, \ldots, n, ( i_2, \ldots, i_k
) \not = ( i, \ldots, i ) \right\}, $$ for $i = 1, \ldots, n$.
\end{theorem}

A substantial portion of this theorem is still true when $\T$ is not
symmetric. As we are only concerned with real symmetric tensors, we
do not go to this in detail.

We call $\sum \{ t_{ii_2 \cdots i_k} : i_2, \ldots, i_k = 1, \ldots,
n, ( i_2, \ldots, i_k ) \not = ( i, \ldots, i ) \}$ the $i$th {\bf
off-diagonal sum} of $\T$.

The set of eigenvalues of $\T$ are called the {\bf spectrum} of
$\T$.  The largest modulus of the eigenvalues of $\T$ is called the
{\bf spectral radius} of $\T$, denoted by $\rho(\T)$.

Following \cite{CPZ}, $\T$ is called {\bf reducible} if there exists
a proper nonempty subset $I$ of $\{ 1, \ldots, n \}$ such that
$$t_{i_1\cdots i_k} = 0, \ \ \forall i_1 \in I, \ \ \forall i_2,
\ldots , i_k \not \in I.$$ If $\T$ is not reducible, then we say
that $\T$ is {\bf irreducible}.    If we take $I = \{ 1, \ldots, n-1
\}$, it is evident that $\LL$ and $\Q$ are reducible.

Suppose that $\T = (t_{i_1 \cdots i_k})$ is a $k$th order
$n$-dimensional tensor.   Construct a graph $\hat G(\T) = (\hat V,
\hat E)$, where $\hat V = \cup_{j=1}^n V_j, V_j$ is a copy of  $\{
1,\ldots, n \}$, for $j = 1, \ldots, n$. Assume that $i_j \in V_j,
i_l \in V_l, j \not = l$.   The edge $(i_j, i_l) \in \hat E$ if and
only if $t_{i_1 \cdots i_k} \not =  0$ for some $k-2$ indices $\{
i_1, \ldots, i_k \} \setminus \{ i_j, i_l \}$.   The tensor $\T$ is
called {\bf weakly irreducible} if $\hat G(\T)$ is connected.   The
original definition in \cite{FGH} for weakly irreducible tensors are
only for nonnegative tensors.   Here we remove the nonnegativity
restriction. As observed in \cite{FGH}, an irreducible tensor is
always weakly irreducible.   Very recently, Pearson and Zhang
\cite{PZ} proved that the adjacency tensor $\A$ is weakly
irreducible if and only if the $k$-graph $G$ is connected. Clearly,
if the adjacency tensor $\A$ is weakly irreducible, then $\LL$ and
$\Q$ are weakly irreducible. This shows that if $G$ is connected,
then $\A, \LL$ and $Q$ are weakly irreducible.

If the entries $t_{i_1\cdots i_k}$ are nonnegative, $\T$ is called a
{\bf nonnegative tensor}.   There is a rich theory on eigenvalues of
a nonnegative tensor \cite{CPZ, CPZ1, FGH, LZI, NQZ, YY, ZQX}. We
now summarize the Perron-Frobenius theorem for nonnegative tensors,
established in \cite{CPZ, FGH, YY}.  With the new definitions of
H$^+$-eigenvalues and H$^{++}$-eigenvalues, this theorem can be
stated concisely.

\begin{theorem} {\bf (The Perron-Frobenius Theorem for Nonnegative
Tensors)} \label{t2}

1. {\bf (Yang and Yang 2010)} If $\T$ is a nonnegative tensor of
order $k$ and dimension $n$, then $\rho(\T)$ is an H$^+$-eigenvalue
of $\T$.

2. {\bf (Friedland, Gaubert and Han 2011)}  If furthermore $\T$ is
weakly irreducible, then $\rho(\T)$ is the unique
H$^{++}$-eigenvalue of $\T$, with the unique eigenvector $x \in
\Re^n_{++}$, up to a positive scaling coefficient.

3. {\bf (Chang, Pearson and Zhang 2008)} If moreover $\T$ is
irreducible, then $\rho(\T)$ is the unique H$^+$-eigenvalue of $\T$.

\end{theorem}

The tensors $\LL$ and $\Q$ are reducible. This permits the
possibility that they have some strict H$^{+}$ eigenvalues. In the
next five sections, we will study their H$^{+}$ eigenvalues.

\section{H$^+$-Eigenvalues of $\A$, $\LL$ and
$Q$}

Theorem \ref{t1} establishes some basic properties of eigenvalues of
the adjacency tensor $\A$, the Laplacian tensor $\LL$ and the
signless Laplacian tensors $\Q$. Note that they are all real $k$th
order $n$-dimensional symmetric tensors.  Both $\A$ and $\Q$ are
nonnegative tensors.   The diagonal entries of $\A$ are zero. The
$i$th diagonal entry of $\LL$ and $\Q$ is $d_i > 0$. All the
off-diagonal entries of $\A$ and $\Q$ are nonnegative. All the
off-diagonal entries of $\LL$ are non-positive.   The $i$th
off-diagonal sum of $\A$ and $\Q$ is $d_i$.  The $i$th off-diagonal
sum of $\LL$ is $-d_i$.

\begin{theorem} {\bf (Basic Properties of Eigenvalues of $\A, \LL$ and $\Q$ )}  \label{t3}

Assume that $k \ge 3$.   The following conclusions hold for
eigenvalues of $\A, \LL$ and $\Q$.

(a). A number $\lambda \in C$ is an eigenvalue of $\A$
(respectively, $\LL$ or $\Q$) if and only if it is a root of the
characteristic polynomial $\phi_\A$ (respectively, $\phi_\LL$ or
$\phi_\Q$).

(b). The number of eigenvalues of $\A$ (respectively, $\LL$ or $\Q$)
is $n(k-1)^{n-1}$. Their product is equal to det$(\A)$
(respectively, det$(\LL)$ or det$(\Q)$).

(c). The sum of all the eigenvalues of $\A$ is zero.  The sum of all
the eigenvalues of $\LL$ or $\Q$ is $(k-1)^{n-1}\sum_{i=1}^n d_i =
k(k-1)^{n-1}m$.

(d). The eigenvalues of $\A$ lie in the disk $\{ \lambda : |\lambda|
\le \Delta \}$. The eigenvalues of $\LL$ and $\Q$ lie in the disk
$\{ \lambda : |\lambda - \Delta| \le \Delta \}$.

(e). $\LL$ and $\Q$ are positive semi-definite when $k$ is even.

\end{theorem}

\noindent {\bf Proof.} The conclusions (a), (b), (c) and (d) follow
directly from Theorem \ref{t1} (a), (b), (c) and (e), and the basic
structure of $\A$. $\LL$ and $\Q$.  By (d), the real parts of all
the eigenvalues of $\LL$ and $\Q$ are nonnegative.  Then (e) follows
from Theorem \ref{t1} (d). \ep

We now discuss H$^+$-eigenvalues of $\LL$.

\begin{theorem} {\bf (H$^+$-Eigenvalues of $\LL$ )}  \label{t4}
Assume that $k \ge 3$. For $j = 1, \cdots, n$, $d_j$ is a strict
H$^+$-eigenvalue of $\LL$ with H-eigenvector $\jj$.   Zero is the
unique H$^{++}$-eigenvalue of $\LL$ with H-eigenvector $\1$, and is
the smallest H-eigenvalue of $\LL$.
\end{theorem}

\noindent {\bf Proof.} A real number $\mu$ is an H-eigenvalue of
$\LL$, with H-eigenvector $x$, if and only if $x \in \Re^n$, $x \not
= 0$, and $\LL x^{k-1} = \mu x^{[k-1]}$, i.e.,
\begin{equation} \label{e3}
d_ix_i^{k-1} - \sum \left\{ {1 \over (k-1)!}x_{i_2}\cdots x_{i_k} :
(i, i_2, \cdots, i_k ) \in E \right\} = \mu x_i^{k-1},
\end{equation}
for $i = 1 \cdots, n$.   We now may easily verify that for $j = 1,
\cdots, n$, $d_j$ is an H$^+$-eigenvalue of $\LL$ with H-eigenvector
$\jj$, and zero is an H$^{++}$-eigenvalue of $\LL$ with
H-eigenvector $\1$.   By Theorem \ref{t3} (d), the real parts of all
the eigenvalues of $\LL$ are nonnegative. Thus, zero is the smallest
H-eigenvalue of $\LL$.  Assume that $x$ is a positive H-eigenvector
of $\LL$, associated with an H-eigenvalue $\mu$. By Theorem \ref{t3}
(d), $\mu \ge 0$. Let $x_j = \min_i \{ x_i \}$. By (\ref{e3}), we
have
$$\mu = d_j - \sum \left\{ {1 \over (k-1)!}{x_{i_2} \over x_j} \cdots {x_{i_k} \over x_j}:
(j, i_2, \cdots, i_k ) \in E \right\} \le d_j - d_j = 0.$$    This
shows that $\mu = 0$.  Thus, zero is the unique H$^{++}$ eigenvalue
of $\LL$, and $d_j$ is a strict H$^+$ eigenvalue of $\LL$, for $j =
1, \cdots, n$. \ep

As in spectral graph theory \cite{BH}, we may call eigenvalues
(respectively, H-eigenvalue or H$^+$-eigenvalue or
H$^{++}$-eigenvalue or spectrum or spectral radius) of $\A$ as
eigenvalues (respectively, H-eigenvalue or H$^+$-eigenvalue or
H$^{++}$-eigenvalue or spectrum or spectral radius) of the $k$-graph
$G$, or simply eigenvalues (respectively, H-eigenvalue or
H$^+$-eigenvalue or H$^{++}$-eigenvalue or spectrum or spectral
radius) if the context is clear. Similarly, we may call eigenvalues
(respectively, H-eigenvalues or H$^+$-eigenvalue or
H$^{++}$-eigenvalue or spectrum or spectral radius) of $\LL$ and $Q$
as Laplacian and signless Laplacian eigenvalues (respectively,
H-eigenvalues or H$^+$-eigenvalue or H$^{++}$-eigenvalue or spectrum
or spectral radius) of $G$, or simply Laplacian and signless
Laplacian eigenvalues (respectively, H-eigenvalues or
H$^+$-eigenvalue or H$^{++}$-eigenvalue or spectrum or spectral
radius) if the context is clear.

Theorem 3.1 of \cite{CD} concerns the spectrum of the union of two
disjoint hypergraphs.   Checking its proof, it also holds for
Laplacian and signless Laplacian spectra.   This will be useful for
our further discussion.  We state it here but omit its proof as the
proof is the same as the proof of Theorem 3.1 of \cite{CD}.

\begin{theorem}  {\bf (The Union of Two Disjoint Hypergraphs)}
\label{t5} Suppose $G=(V, E)$ is the union of two disjoint
hypergraphs $G_1 = (V_1, E_1)$ and $G_2 = (V_2, E_2)$, where $|V_1|
= n_1, |V_2| = n_2, n_1 + n_2 = n = | V |$.   Then the spectrum
(respectively, the Laplacian spectrum or the signless Laplacian
spectrum) of $G$ is the union of the spectra (respectively, the
Laplacian spectra or the signless Laplacian spectra) of $G_1$ and
$G_2$, where, as multisets, an eigenvalue with multiplicity $r$ in
the spectrum (respectively, the Laplacian spectrum or the signless
Laplacian spectrum) of $G_1$ occurs in the spectrum (respectively,
the Laplacian spectrum or the signless Laplacian spectrum) of $G$
with multiplicity $r(k-1)^{n_2}$.
\end{theorem}

In general, $G$ may be decomposed into components $G_r = (V_r, E_r)$
for $r = 1, \ldots , s$.   If $s = 1$, then $G$ is connected. Denote
the adjacency tensor and the signless Laplacian tensor of $G_r$ by
$\A(G_r)$ and $\Q(G_r)$ respectively, for $r = 1, \ldots, s$. Then
by Theorem \ref{t5},
$$\rho(\A) = \max_{r=1,\cdots, s} \{ \rho(\A(G_r)) \},\ \  \rho(\Q) = \max_{r=1,\ldots, s} \{ \rho(\Q(G_r)) \}.$$
With the above discussion, we are now ready to study
H$^+$-eigenvalues of $\Q$ and $\A$.

\begin{theorem} {\bf (H$^+$-Eigenvalues of $\Q$ )}  \label{t6}
Assume that $k \ge 3$.   Suppose that $G$ has $s$ components $G_r =
(V_r, E_r)$ for $r = 1, \ldots , s$.  For $j = 1, \ldots, n$, $d_j$
is a strict H$^+$-eigenvalue of $\Q$ with an H-eigenvector $\jj$.
Let $\nu_1 = \rho(\Q)$.  If $\nu_1 \equiv \rho(\Q(G_r))$ for $r = 1,
\ldots, s$, then $\nu_1$ is the unique H$^{++}$-eigenvalue of $\Q$.
Otherwise, $\Q$ has no H$^{++}$-eigenvalue, and for $r = 1, \ldots,
s$, $\rho(\Q(G_r))$ is a strict H$^+$-eigenvalue of $\Q$.
\end{theorem}

\noindent {\bf Proof.} A real number $\nu$ is an H-eigenvalue of
$\Q$, with an H-eigenvector $x$, if and only if $x \in \Re^n$, $x
\not = 0$, and $\Q x^{k-1} = \nu x^{[k-1]}$, i.e.,
\begin{equation} \label{e4}
d_ix_i^{k-1} + \sum \left\{ {1 \over (k-1)!}x_{i_2}\cdots x_{i_k} :
(i, i_2, \ldots, i_k ) \in E \right\} = \nu x_i^{k-1},
\end{equation}
for $i = 1, \ldots, n$.   Then, we may easily verify that for $j =
1, \ldots, n$, $d_j$ is an H$^+$-eigenvalue of $Q$ with an
H-eigenvector $\jj$.

For $r = 1, \ldots, s$, as $G_r$ is connected, $\Q(G_r)$ is weakly
irreducible by \cite{PZ}.   By Theorem \ref{t2}, $\rho(\Q(G_r))$ is
the unique H$^{++}$-eigenvalue of $\Q(G_r)$, with a positive
H-eigenvector $x^{(r)} \in \Re^{|V_r|}$.   In (\ref{e4}), let $\nu =
\rho(\Q(G_r))$, $x_i = x^{(r)}_i$ if $i \in V_r$ and $x_i = 0$ if $i
\not \in V_r$.  Then we see that (\ref{e4}) is satisfied for $i = 1,
\ldots, n$.   This shows that for $r = 1, \ldots, s$,
$\rho(\Q(G_r))$ is an H$^+$-eigenvalue of $\Q$.

Assume that $\nu$ is an H$^{++}$-eigenvalue of $\Q$ with a positive
H-eigenvector $x$.  For $r = 1, \ldots, s$, define $x^{(r)} \in
\Re^{|V_r|}$ by $x^{(r)}_i = x_i$ for $i \in V_r$.  Then $x^{(r)}$
is a positive H-eigenvector in $\Re^{|V_r|}$.   By (\ref{e4}), $\nu$
is an H$^{++}$-eigenvalue of $\Q(G_r)$.   Since $\Q(G_r)$ is weakly
irreducible, by Theorem \ref{t2}, $\nu = \rho(\Q(G_r))$.  Thus, if
$\Q$ has an H$^{++}$-eigenvalue, then it must be $\nu_1 = \rho(\Q)
\equiv \rho(\Q(G_r))$ for $r = 1, \ldots, s$.   This completes our
proof.  \ep

\begin{theorem} {\bf (H$^+$-Eigenvalues of $\A$ )}  \label{t7}
Assume that $k \ge 3$.   Then zero is a strict H$^+$-eigenvalue of
$\A$. Suppose that $G$ has $s$ components $G_r = (V_r, E_r)$ for $r
= 1, \ldots , s$.  Let $\lambda_1 = \rho(\A)$.  If $\lambda_1 \equiv
\rho(\A(G_r))$ for $r = 1, \ldots, s$, then $\lambda_1$ is the
unique H$^{++}$-eigenvalue of $\A$. Otherwise, $\A$ has no
H$^{++}$-eigenvalue, and for $r = 1, \ldots, s$, $\rho(\A(G_r))$ is
a strict H$^+$-eigenvalue of $\A$.
\end{theorem}

\noindent {\bf Proof.} Zero is an H-eigenvalue of $\A$, with an
H-eigenvector $x$, if and only if $x \in \Re^n$, $x \not = 0$, and
$\A x^{k-1} = 0$, i.e.,
$$\sum \left\{ {1 \over (k-1)!}x_{i_2}\cdots x_{i_k} : (i,
i_2, \ldots, i_k ) \in E \right\} = 0,$$ for $i = 1 \cdots, n$. Let
$x$ be a vector in $\Re^n_+$ with $1 \le {\rm supp}(x) \le k-2$.
Then we see that $x$ is a nonnegative H-eigenvector of $\A$,
corresponding to the zero H-eigenvalue.  Thus, zero is an
H$^+$-eigenvalue of $\A$.  The proof of the remaining conclusions of
this theorem is similar to the last part of the proof of the last
theorem.   We omit it.  \ep

For some $k$-graph $G$, $\LL, \Q$ and $\A$ may have more strict
H$^+$-eigenvalues.   For example, let $k=3, n=8, m=8$, and $E = \{
(1, 2, 3), (1, 4, 5), (2, 4, 5), (3, 4, 5), (4, 5, 6), (4, 5, 7),
(4, 5, 8), (6, 7, 8) \}$.   Then $d_1 = d_2 = d_3 = d_6 = d_7 = d_8
= 2$ and $d_4 = d_5 = 6$ are strict H$^+$-eigenvalues of $\LL$ and
$\Q$, $0$ is a strict H$^+$-eigenvalue of $\A$.   It is easy to
verify that $\mu = 1$, $\nu = 3$ and $\lambda = 1$ are also strict
H$^+$-eigenvalues of $\LL, \Q$ and $\A$, with an H-eigenvector $(1,
1, 1, 0, 0, 0, 0, 0)$.

We will not identify all strict H$^+$-eigenvalues of $\LL, \Q$ and
$\A$, but we will identify the largest and the smallest
H$^+$-eigenvalues of $\LL$ and $\Q$, and establish their maximum or
minimum properties in the next few sections.  They are the most
important H$^+$-eigenvalues of $\LL$ and $\Q$.

There are also H-eigenvalues of $\LL$ and  $\Q$ which are not
H$^+$-eigenvalues.   We will give such an example in Sections 5 and
7.

Theorems \ref{t4}, \ref{t6} and \ref{t7} say that each of $\LL, \Q$
and $\A$ has at most one H$^{++}$-eigenvalue.   Actually, a real
symmetric matrix has at most one H$^{++}$-eigenvalue.  By Theorem
\ref{t2}, a weakly irreducible nonnegative tensor has at most one
H$^{++}$-eigenvalue.   By extending the proof of Theorem \ref{t6},
probably this is also true for a general nonnegative tensor.  We may
also show that this is true for a real diagonal tensor. However, by
numerical experiments, we found that this is not true for some real
symmetric tensors.   Thus, we ask the following question.

{\bf Question 1.} Is there a reasonable class of real symmetric
tensors, which includes the above cases, such that any tensor in
this class has at most one H$^{++}$-eigenvalue?

%For the other complex eigenvalues of $\A, \LL$ and $\Q$, we may
%define an order for them.   Let complex numbers $\lambda_i =
%\alpha_i + \beta_i\sqrt{-1}$, where $\alpha_i$ and $\beta_i$ are
%real numbers, for $i = 1, 2$.   If $\alpha_1 > \alpha_2$, or
%$\alpha_1 = \alpha_2$ and $\beta_1 > \beta_2$, then we denote
%$\lambda_1 \prec \lambda_2$ or $\lambda_2 \succ \lambda_1$. If
%$\lambda_1 \prec \lambda_2$ or $\lambda_1 = \lambda_2$, then we
%denote $\lambda_1 \preceq \lambda_2$ or $\lambda_2 \succeq
%\lambda_1$.    Then we may order the complex eigenvalues of $\A,
%\LL$ and $\Q$ as $\rho(\A) = \lambda_1 \succeq \lambda_2 \succeq
%\cdots \succeq \lambda_b$, $0 = \mu_1 \preceq \mu_2 \preceq \cdots
%\preceq \mu_b$, and  $\rho(\Q) = \nu_1 \succeq \nu_2 \succeq \cdots
%\succeq \nu_b$, where $b = n(k-1)^{n-1}$.

%In \cite{CD}, it was proved that $\bar d \le \lambda_1 \le \Delta$.
%If $G$ is a $d$-regular $k$-graph, we have $\lambda_1 = d$,
%$\lambda_i = d - \mu_i = \nu_i - d$ for $i = 1, \cdots, b$.   This
%is similar to the results for graphs in Section 1.3.1 of \cite{BH}.

\section{The Smallest Laplacian H-Eigenvalue}

The smallest Laplacian H-eigenvalue of $G$ is $\mu_1 = 0$.   By
Theorem \ref{t4}, $\1$ is an H-eigenvector of $\LL$, associated with
the H$^{++}$-eigenvalue $\mu_1 = 0$.  We say that $x \in \Re^n$ is a
{\bf binary vector} if $x_i$ is either $0$ or $1$ for $i = 1,
\ldots, n$. Thus, $\1$ is a binary H-eigenvector of $\LL$,
associated with the H-eigenvalue $\mu_1 = 0$.   We say that a binary
H-eigenvector $x$ of $\LL$, associated with an H-eigenvalue $\mu$,
is a {\bf minimal binary H-eigenvector} of $\LL$, associated with
$\mu$, if there does not exist another binary H-eigenvector $y$ of
$\LL$, associated with $\mu$, such that supp($y$) is a proper subset
of supp($x$).

Let $G = (V, E)$ be a $k$-graph. For $e_p = (i_1, \ldots, i_k) \in
E$, define a $k$th order $n$-dimensional symmetric tensor $\LL(e_p)$
by
$$\LL(e_p)x^k = \sum_{j=1}^k x_{i_j}^k - k x_{i_1}\cdots x_{i_k}$$
for any $x \in C^n$.   Then, for any $x \in C^n$, we have
$$\LL x^k = \sum_{e_p \in E} \LL(e_p)x^k.$$

\begin{Theorem}  {\bf (The Smallest Laplacian H-Eigenvalue)} \label{t8}
For a $k$-graph $G$, we have the following conclusions.

(a). For any $x \in \Re^+$, $\LL x^k \ge 0$.   We have $$0 = \min \{
\LL x^k : x \in \Re^n_+, \sum_{i=1}^n x_i^k = 1 \}.$$

(b). A binary vector $x \in \Re^n$ is a minimal binary H-eigenvector
of $\LL$ associated with the H-eigenvalue $\mu_1 = 0$ if and only if
supp($x$) is the vertex set of a component of $G$.

(c).  A vector $x \in \Re^n$ is an H-eigenvector of $\LL$ associated
with the H-eigenvalue $\mu_1 = 0$ if it is a nonzero linear
combination of minimal binary H-eigenvectors of $\LL$ associated
with the H-eigenvalue $\mu_1 = 0$.

\end{Theorem}

\noindent {\bf Proof.} (a). For any $e_p = (i_1, \ldots, i_k) \in E$
and $x \in \Re^n_+$, we know that the arithmetic mean of $x_{i_1}^k,
\cdots, x_{i_k}^k$ is greater than or equal to their geometric mean,
i.e.,
$${1 \over k}\sum_{j=1}^k x_{i_j}^k \ge x_{i_1}\cdots x_{i_k}.$$
This implies that $\LL (e_p)x^k \ge 0$.   Thus, $\LL x^k \ge 0$ for
any $x \in \Re^n_+$.   As $\LL y^k = 0$, where $y = {l \over n^{1
\over k}}$, we have
$$0 = \min \{ \LL x^k : x \in \Re^n_+, \sum_{i=1}^n x_i^k = 1 \}.$$

(b). A nonzero vector $x \in \Re^n$ is an H-eigenvector of $\LL$,
associated the the H-eigenvalue $\mu_1 = 0$, if and only $\LL
x^{k-1} =0$, i.e.,
\begin{equation} \label{mu1}
d_ix_i^{k-1} = \sum \left\{ {1 \over (k-1)!}x_{i_2}\cdots x_{i_k} :
(i, i_2, \ldots, i_k) \in E \right\},
\end{equation}
for $i = 1, \ldots, n$.

Suppose that $x$ is a binary vector and supp($x$) is the vertex set
of a component of $G$.   Then the equation (\ref{mu1}) reduces to
$d_i = d_i$ if $i \in$ supp($x$), and $0 = 0$ if $i \not \in$
supp($x$). Thus, $x$ is a binary H-eigenvector of $\LL$ associated
with H-eigenvalue $\mu_1 = 0$. Suppose that $y$ is a binary vector
and supp($y$) is a proper subset of supp($x$).  Then there are $i
\in$ supp($y$) and an edge $(i, i_2, \ldots, i_k) \in E$ such that
one of the indices $\{ i_2, \ldots, i_k \}$ not in supp($y$). Then,
for this $i$, by replacing $x$ by $y$ in (\ref{mu1}), the left hand
side of (\ref{mu1}) becomes $d_i$, while the right hand side of
(\ref{mu1}) is strictly less than $d_i$, i.e., (\ref{mu1}) does not
hold under this replacement.  This shows that $y$ cannot be a binary
H-eigenvector of $\LL$ associated with H-eigenvalue $\mu_1 = 0$,
i.e., $x$ is a minimal binary H-eigenvector of $\LL$ associated with
the H-eigenvalue $\mu_1 = 0$.

On the other hand, suppose that $x$ is a binary H-eigenvector of
$\LL$ associated with the H-eigenvalue $\mu_1 = 0$.  Let $i \in$
supp($x$).  Then, in order that the equation (\ref{mu1}) holds for
$i$, for any $(i, i_2, \ldots, i_k) \in E$, we must have $i_2,
\ldots, i_k \in$ supp($x$).   This shows that supp($x$) is either
the vertex set of a component of $G$, or the union of the vertex
sets of several components of $G$.   This proves (b).

(c). Let $\left\{ y^{(1)}, \ldots, y^{(s)} \right\}$ be the set of
binary H-eigenvectors of $\LL$ associated with H-eigenvalue $\mu_1 =
0$.

Suppose that $x$ is a nonzero linear combination of  $y^{(1)},
\ldots, y^{(s)}$, $x = \sum_{r=1}^s \alpha_r y{(r)}$, where
$\alpha_r$ are real numbers.    If $i \in$ supp($y^{(r)}$) for some
$r$, then the equation (\ref{mu1}) is $\alpha_r^{k-1}d_i =
\alpha_r^{k-1}d_i$.  Otherwise, the equation (\ref{mu1}) is $0=0$.
Thus, $x$ is an H-eigenvector of $\LL$ associated with the
H-eigenvalue $\mu_1 = 0$. This proves (c). \ep

%As in \cite{HQ}, we may call the dimension of $V$ the {\bf
%geometrical multiplicity} of $\mu_1 = 0$.  Suppose that the
%eigenvalues of $\LL$ are $0 = \mu_1 = \cdots = \mu_s \prec \mu_{s+1}
%\preceq \dots \preceq \mu_b$, where $b = n(k-1)^{n-1}$.  Then $s$ is
%called the {\bf algebraic multiplicity} of $\mu_1=0$.  By Theorem
%\ref{t5}, the geometrical multiplicity of $\mu_1 = 0$ is not greater
%than $n \over k$.    By Theorem \ref{t4}, the algebraic multiplicity
%of $\mu_1$ can be much bigger.  Thus, they are not the same in
%general.   The question is if the geometrical multiplicity of $\mu_1
%= 0$ is one, will the algebraic multiplicity of $\mu_1$ also be one?

\begin{Corollary}\label{c1} The following two statements are equivalent.

(a). The $k$-graph $G$ is connected.

(b). The vector $\1$ is the unique minimal binary H-eigenvector of
$\LL$ associated with the H-eigenvalue $\mu_1 = 0$.

\end{Corollary}

\section{The Largest Laplacian H$^+$-Eigenvalue}

In Section 3, we showed that zero is the unique Laplacian
H$^{++}$-eigenvalue of $G$, and $d_j$ is a strict H$^+$-eigenvalue
of $G$, for $j = 1, \ldots, n$.   We now identify the largest
Laplacian H$^+$-eigenvalue of $G$, and establish a maximum property
of this Laplacian H$^+$-eigenvalue.

\begin{theorem} {\bf (The Largest Laplacian H$^+$-Eigenvalue)}
\label{t9} Assume that $k \ge 3$.   The largest Laplacian
H$^+$-eigenvalue of $G$ is $\Delta = \max_i \{ d_i \}$. We have
\begin{equation} \label{e5}
\Delta = \max \{ \LL x^k : x \in \Re^n_+, \sum_{i=1}^n x_i^k = 1 \}.
\end{equation}
\end{theorem}

\noindent {\bf Proof.} Suppose that $\mu$ is a Laplacian
H$^+$-eigenvalue of $G$ associated with nonnegative H-eigenvector
$x$.   Assume that $x_j > 0$. By (\ref{e3}), we have
$$\mu x_j^{k-1}
= d_jx_j^{k-1} - \sum \left\{ {1 \over (k-1)!}x_{i_2}\cdots x_{i_k}
: (i, i_2, \ldots, i_k ) \in E \right\} \le d_jx_j^{k-1}.$$ This
implies that
$$\mu \le d_j \le \Delta.$$
By Theorem \ref{t4}, $\Delta$ is an H$^{+}$-eigenvalue of $\LL$.
Thus, $\Delta$ is the largest H$^{+}$-eigenvalue of $\LL$.

Suppose that $\Delta = d_j$.   Let $x = \jj$.   Then $x$ is a
feasible point of the maximization problem in (\ref{e5}).   We have
$$\LL x^k = \sum_{i = 1}^n \left[ d_ix_i^k - \sum \left\{ {1 \over (k-1)!}x_ix_{i_2}\cdots x_{i_k}
: (i, i_2, \ldots, i_k ) \in E \right\} \right] = \Delta.$$ This
shows that
$$\Delta \le \max \{ \LL x^k : x \in \Re^n_+, \sum_{i=1}^n x_i^k = 1
\}.$$

On the other hand, suppose $x^*$ is a maximizer of the maximization
problem in (\ref{e5}).   As the feasible set is compact, and the
objective function is continuous, such a maximizer exists.  By
optimization theory, for $i = 1, \cdots, n$, either $x_i^* = 0$ and
\begin{equation} \label{e6}
d_i(x_i^*)^{k-1} - \sum \left\{ {1 \over (k-1)!}x^*_{i_2}\cdots
x^*_{i_k} : (i, i_2, \ldots, i_k ) \in E \right\} \ge \mu
(x_i^*)^{k-1},
\end{equation}
or $x_i^* > 0$ and
\begin{equation} \label{e7}
d_i(x_i^*)^{k-1} - \sum \left\{ {1 \over (k-1)!}x^*_{i_2}\cdots
x^*_{i_k} : (i, i_2, \ldots, i_k ) \in E \right\} = \mu
(x_i^*)^{k-1},
\end{equation} where $\mu$ is a Lagrange multiplier.  As $x^*$ is feasible for the maximization problem,
(\ref{e7}) holds for at least one $i$, say $i_0$.  We have
$$d_{i_0}(x_{i_0}^*)^{k-1} \ge \mu (x_{i_0}^*)^{k-1}.$$
As $x_{i_0}^* > 0$, we have $\mu \le d_{i_0} \le \Delta$.
Multiplying (\ref{e6}) and (\ref{e7}) by $x_i^*$ and summing up them
for $i = 1, \ldots , n$, we have
$$\LL (x^*)^k = \mu \sum_{i=1}^n (x^*_i)^k = \mu.$$
Thus,
$$\mu = \max \{ \LL x^k : x \in \Re^n_+, \sum_{i=1}^n x_i^k = 1
\}.$$ This shows that
$$\Delta \ge \max \{ \LL x^k : x \in \Re^n_+, \sum_{i=1}^n x_i^k = 1
\}.$$ Hence, (\ref{e5}) holds. \ep

In general, $\Delta$ may not be the largest H-eigenvalue of $\LL$.
For example, let $n=k=6, m=1$ and $E= \{ (1, 2, 3, 4, 5, 6) \}$.
Then $\Delta = 1$, while $\mu = 2$ is an H-eigenvalue of $\LL$ with
an H-eigenvector $(1, 1, 1, -1, -1, -1)$.

\section{The Largest H-Eigenvalue and The Largest Signless Laplacian H-Eigenvalue}

The largest H-eigenvalue is $\lambda_1 = \rho(\A)$.  The largest
signless Laplacian H-eigenvalue is $\nu_1 = \rho(\Q)$.    As both
$\A$ and $\Q$ are nonnegative tensors, their properties are similar.
We thus discuss them together.

When $k$ is even, by \cite{Qi}, we know that
$$\lambda_1 = \max \{ \A x^k : x \in \Re^n, \sum_{i=1}^n x_i^k = 1
\},$$ and
$$\nu_1 = \max \{ \Q x^k : x \in \Re^n, \sum_{i=1}^n x_i^k = 1
\}.$$ The feasible sets of the above two maximization problems are
the same.   It is a compact set when $k$ is even.   When $k$ is odd,
it is not compact.   We intend to establish some maximum properties
of $\lambda_1$ and $\nu_1$, which hold whenever $k$ is even or odd.

Corollary 3.4 of \cite{CD} indicates that when $G$ is connected,
\begin{equation} \label{e6.1}
\lambda_1 = \max \{ \A x^k : x \in \Re^n_+, \sum_{i=1}^n x_i^k = 1
\}.
\end{equation}
Using a similar argument, we may show that when $G$ is connected,
\begin{equation} \label{e6.2}
\nu_1 = \max \{ \Q x^k : x \in \Re^n_+, \sum_{i=1}^n x_i^k = 1 \}.
\end{equation}
We wish to show that (\ref{e6.1}) and (\ref{e6.2}) hold even if $G$
is not connected.

\begin{theorem} {\bf (The Largest H-Eigenvalue and the Largest Signless Laplacian H-Eigenvalue)}
\label{t10} Assume that $k \ge 3$.  Then (\ref{e6.1}) and
(\ref{e6.2}) always hold.
\end{theorem}

\noindent {\bf Proof.} We now prove (\ref{e6.1}).   Suppose that $G$
is decomposed to some components $G_r = (V_r, E_r)$ for $r = 1,
\ldots, s$.   Then $\lambda_1 = \max \{ \rho(\A(G_r)) : r = 1,
\ldots, s \}$, and for $r = 1, \ldots, s$,
$$\rho(\A(G_r)) = \max \{ \A(G_r) (x^{(r)})^k : x^{(r)} \in \Re^{|V_r|}_+, \sum_{i\in V_r} \left(x_i^{(r)}\right)^k = 1
\}.$$ Suppose that $\lambda_1 = \rho(\A(G_j))$ for some $j$. Define
$x \in \Re^n_+$ by $x_i = x^{(r)}_i$ if $i \in V_r$ and $x_i = 0$
otherwise.   Then $\sum_{i=1}^n x_i^k = 1$, and $\A x^k = \A(G_j)
(x^{(j)})^k$.   We see that $\lambda_1 = \A x^k$ and $x$ is a
feasible point of the maximization problem in (\ref{e6.1}).  This
shows that
$$\lambda_1 \le \max \{ \A x^k : x \in \Re^n_+,
\sum_{i=1}^n x_i^k = 1 \}.$$

On the other hand, suppose that $x_*$ is a maximizer of the
maximization problem in (\ref{e6.1}).   Then,
$$\max \{ \A x^k : x \in \Re^n_+,
\sum_{i=1}^n x_i^k = 1 \} = \A x_*^k = \sum_{r=1}^s \A(G_r) (\bar
x^{(r)})^k,$$ where $\bar x^{(r)} \in \Re^{|V_r|}_+$ and $\bar
x^{(r)}_i = (x_*)_i$ for $i \in V_r$, for $r = 1, \ldots, s$.   For
$r = 1, \ldots, s$, assume that $\alpha_r = \sum_{i \in V_r}
(x_*)_i^k$.  Then $\alpha_r \ge 0$ for $r = 1, \ldots, s$, and
$\sum_{r=1}^s \alpha_r = 1$.   If $\alpha_r > 0$, then define
$x^{(r)} \in \Re^{|V_r|}_+$ by $x^{(r)} = {1 \over (\alpha_r)^{1
\over k}} \bar x^{(r)}$.   Then $\sum_{i \in V_r}
\left(x_i^{(r)}\right)^k = 1$.   We now have
$$\A x_*^k = \sum \{ \A(G_r) (\bar
x^{(r)})^k : \alpha_r > 0 \}= \sum \{ \alpha_r \A(G_r) (x^{(r)})^k :
\alpha_r > 0 \} \le \sum \{ \alpha_r \rho(\A(G_r)) : \alpha_r > 0
\}$$
$$\le \sum \{ \alpha_r \lambda_1 : \alpha_r > 0
\} = \lambda_1.$$ Thus, we have
$$\lambda_1 \ge \max \{ \A x^k : x \in \Re^n_+,
\sum_{i=1}^n x_i^k = 1 \}.$$ Hence, (\ref{e6.1}) holds.

Similarly, we may show that (\ref{e6.2}) holds. \ep

\begin{Corollary} \label{c2} {\bf (Bounds for $\nu_1$)}   We always have
\begin{equation} \label{e6.3}
\max \{ \Delta, 2\bar d \} \le \nu_1 \le 2\Delta.
\end{equation}
\end{Corollary}

\noindent {\bf Proof.}  By Theorem \ref{t3} (d), we have that
$$0  \le \nu_1 \le 2\Delta.$$
In (\ref{e6.2}), letting $x = {l \over n^{1 \over k}}$, we see that
$\nu_1 \ge 2\bar d$.  Assume that $d_j = \Delta$.  In (\ref{e6.2}),
letting $x = \jj$, we see that $\nu_1 \ge \Delta$.  Thus, we always
have
$$\nu_1 \ge \max \{ \Delta, 2\bar d \}.$$
These prove (\ref{e6.3}). \ep

It was established in \cite{CD} that $\bar d \le \lambda_1 \le
\Delta$.

{\bf Question 2}. Are there any formulas related to $\lambda_1$ and
$\nu_1$?

We may compare $\nu_1$, $\lambda_1$ and $\rho(\LL)$.   We prove a
lemma first.

\begin{Lemma} \label{l1}
 If $\C$ is a nonnegative tensor of order $k$ and dimension $n$,
and $\B$ is a tensor of order $k$ and dimension $n$, satisfying
$|\B| \le \C$, then $\rho(\B) \le \rho(\C)$.
\end{Lemma}

\noindent {\bf Proof.}  Let $\C_\epsilon = \C + \epsilon J$, with
$\epsilon > 0$.  Then $\C_\epsilon$ is a positive tensor, thus
irreducible, and $|\B| \le \C_\epsilon$.  By Lemma 3.2 of \cite{YY},
we have $\rho(\B) \le \rho(\C_\epsilon)$. Let $\epsilon \to 0$. As
the eigenvalues of a tensor are roots of the characteristic
polynomial, whose coefficients are polynomials in the entries of
that tensor \cite{Qi}, the spectral radius of that tensor is
continuous in its entries. Then we have $\rho(\B) \le \rho(\C)$. \ep

With this lemma, we immediately have the following proposition.

\begin{Proposition} \label{p10}
For a $k$-graph $G$, we have
$$\nu_1 = \rho(\Q) \ge
\rho(\LL),\ {\rm and}\  \nu_1 = \rho(\Q) \ge \lambda_1 = \rho(\A).$$
\end{Proposition}

Note that it is possible that $\nu_1 = \rho(\Q) = \rho(\LL)$. For
example, let $n=k=6$, $m=1$ and $E = \{ (1, 2, 3, 4, 5, 6)\}$. Then
$G$ is connected.   Thus, $\A, \LL$ and $\Q$ are weakly irreducible.
We have $\LL x^6 = \sum_{i=1}^6 x_i^6 - 6x_1\cdots x_6$ and $\Q x^6
= \sum_{i=1}^6 x_i^6 - 6x_1\cdots x_6$.  We see that $\nu = 2$ is an
H$^{++}$ eigenvalue of $\Q$ with an H-eigenvector $l = (1, 1, 1, 1,
1, 1)$.  By Theorem \ref{t2} (b), we have $\rho(\Q) = 2$.   On the
other hand, we see $\mu = 2$ is an H-eigenvalue of $\LL$ with an
H-eigenvector $l = (1, 1, 1, -1, -1, -1)$.  By Proposition
\ref{p10}, we have $\rho(\LL) = \rho(\Q) = 2$. Thus, it is a
research topic to identify the conditions under which $\rho(\LL) =
\rho(\Q)$.

\smallskip

We now discuss algorithms for computing $\nu_1$.   As $\Q$ is a
nonnegative tensor, we may use algorithms for finding the largest
eigenvalue of a nonnegative tensor to compute it.   However, the
convergence of the NQZ algorithm \cite{NQZ} needs the condition that
$\Q$ is primitive \cite{CPZ1}, and the convergence of the LZI
algorithm needs the condition that $\Q$ is irreducible \cite{LZI}.
These conditions are somewhat strong. The linear convergence of the
LZI algorithm needs the condition that $\Q$ is weakly positive
\cite{ZQX}. A nonnegative tensor $\T = (t_{i_1\cdots i_k})$ is
weakly positive if $t_{ij\cdots j} > 0$ for all $i \not = j, i, j =
1, \ldots, n$. We see that $\Q$ cannot be weakly positive. Thus, it
may not be a good choice to use these two algorithms for computing
$\nu_1$.  Instead, one may use the HHQ algorithm proposed in
\cite{HHQ} to compute $\nu_1$. The HHQ algorithm is globally
R-linearly convergent if $Q$ is weakly irreducible in the sense of
\cite{FGH}. As discussed above, if $G$ is connected, then $Q$ is
weakly irreducible. Thus, the HHQ algorithm is practicable for
computing $\nu_1$ when $G$ is connected.  If $G$ is not connected,
the HHQ algorithm may be used for components (and then the maximum
value chosen), by the observation at the beginning of the proof of
Theorem 11.   This argument is also valid for computing $\lambda_1$.

%\noindent {\bf Proof.} Suppose that $G$ is connected.  Clearly, if $\A$ is weakly irreducible, then $\Q$ is weakly irreducible.
% Thus, it suffices to prove that $\A$ is weakly irreducible.
%Construct $\hat G(\A)$
% and still use the above symbols $V_j, V_i, i_j, i_l$, etc.   Suppose that $j \not = l, i_j \not = i_l$.
%  As $G$ is connected, it is easy to see that the vertices $i_j \in V_j$ and $i_l \in V_l$ are connected.   Now, via a third vertex,
%  we may easily see that any two vertices $i_j, i_l \in V_j, i_j \not = i_l$ are connected, and any two vertices $i_j \in V_j, i_l \in V_l, j \not = l$,
%  but $i_j = i_l$,  are also  connected.   This shows that $\hat G(\A)$ is connected, i.e., $\A$ is weakly irreducible.  \ep
Thus, we may use the HHQ algorithm to compute $\lambda_1$ and
$\nu_1$, and we have global R-linear convergence.

\section{The Smallest Signless Laplacian H$^+$-Eigenvalue}

We now identify the smallest signless Laplacian H$^+$-eigenvalue of
$G$, and establish a minimum property of this signless Laplacian
H$^+$-eigenvalue.

\begin{Theorem} \label{t11} {\bf (The Smallest Signless Laplacian H$^+$-Eigenvalue)}  The smallest signless Laplacian H$^+$-eigenvalue of
$G$ is $\delta$.   We always have
\begin{equation} \label{e7.1}
\delta  =  \min \{ \Q x^k : x \in \Re^n_+, \sum_{i=1}^n x_i^k = 1
\}.
\end{equation}
\end{Theorem}

\noindent {\bf Proof.}   Suppose that $\nu$ is an H$^+$-eigenvalue
of $\Q$, with a nonnegative H-eigenvector $x$.  Suppose that $x_j >
0$.  By (\ref{e4}), we have
$$d_jx_j^{k-1} + \sum \left\{ {1 \over (k-1)!}x_{i_2}\cdots x_{i_k} :
(j, i_2, \ldots, i_k ) \in E \right\} = \nu x_j^{k-1}.$$ This
implies that $d_jx_j^{k-1} \le \nu x_j^{k-1}$, i.e., $\nu \ge d_j
\ge \delta$.   As $\delta$ is an H$^+$-eigenvalue of $\Q$ by Theorem
\ref{t6}, this shows that $\delta$ is the smallest H$^+$-eigenvalue
of $\Q$.

We now prove (\ref{e7.1}).  Suppose that $d_j = \delta$.  Let $x =
\jj$ in (\ref{e7.1}). Then we have
\begin{equation} \label{e7.2}
\delta \ge \min \{ \Q x^k : x \in \Re^n_+, \sum_{i=1}^n x_i^k = 1
\}.
\end{equation}

Suppose that $x^*$ is an optimal solution of the minimization
problem in (\ref{e7.1}). By the optimization theory, there are
Lagrange multipliers $u \in \Re^n$ and $\nu \in \Re$ such that for
$i = 1, \cdots, n$,
\begin{equation} \label{e7.3}
\left(\Q\left(x^*\right)^{k-1}\right)_i = \nu \left(x^*_i\right)^{k-1} + u_i,
\end{equation}
$$x_i^* \ge 0, \ u_i \ge 0, \ x_i^*u_i = 0$$
and
\begin{equation} \label{e7.4}
\sum_{i=1}^n \left(x_i^*\right)^k = 1.
\end{equation}
Let $I =$ supp$(x^*)$.  By (\ref{e7.4}), $I \not = \emptyset$. Then
for $i \in I, u_i = 0$ and for $i \not \in I, x_i^* = 0$.
Multiplying (\ref{e7.3}) by $x_i^*$ and summing from $i=1$ to $n$,
we have
$$\nu = \Q\left(x^*\right)^k.$$
Now assume that $x_j^* = \max \{ x_i^* : i \in I \}$.  Then $x_j^*
> 0$ and $u_j = 0$.  By (\ref{e7.3}), we have
$$\left(\Q\left(x^*\right)^{k-1}\right)_j =
\nu\left(x_j^*\right)^{k-1},$$
which implies that
$$d_j\left(x_j^*\right)^{k-1} \le \nu \left(x_j^*\right)^{k-1}.$$
Thus, $$\nu = \Q \left(x^*\right)^k \ge d_j \ge \delta.$$ Hence,
$$\delta \le \min \{ \Q x^k : x \in \Re^n_+, \sum_{i=1}^n x_i^k =
1 \}.$$ Combining this with (\ref{e7.2}), we have (\ref{e7.1}).
\ep

In general, $\delta$ may not be the smallest H-eigenvalue of $\Q$.
For example, let $n=k=6, m=1$ and $E= \{ (1, 2, 3, 4, 5, 6) \}$.
Then $\delta = 1$, while $\nu = 0$ is an H-eigenvalue of $\Q$ with
an H-eigenvector $(1, 1, 1, -1, -1, -1)$.   In general, we may show
that $\Q$ has a zero H-eigenvalue if and only if $k = 4j+2$ for some
integer $j$, and there is a vector $x \in \Re^k$ such that for any
edge $e_p = (i_1, \ldots, i_k) \in E$, half of $x_{i_1}, \ldots,
x_{i_k}$ are $j$, and the other half are $-1$.  Hence, if $k = 4j$
or if $k = 4j+2$ but such an $x$ does not exist, then $\Q$ is
positive definite.

We now give an application of Theorem \ref{t11}.  Suppose that $S$
is a proper nonempty subset of $V$.   Denote $\bar S = V \setminus
S$. Then $\bar S$ is also a proper nonempty subset of $V$.   The
edge set $E$ is now partitioned into three parts $E(S), E(\bar S)$
and $E(S, \bar S)$.   The edge set $E(S)$ consists of edges whose
vertices are all in $S$.   The edge set $E(\bar S)$ consists of
edges whose vertices are all in $\bar S$.   The edge set $E(S, \bar
S)$ consists of edges whose vertices are in both $S$ and $\bar S$.
We call $E(S, \bar S)$ an {\bf edge cut} of $G$.   If we delete
$E(S, \bar S)$ from $G$, then $G$ is separated into two $k$-graphs
$G[S] = (S, E(S))$ and $G[\bar S] = (\bar S, E(\bar S))$.   For a
vertex $i \in S$, we denote its degree at $G[S]$ by $d_i(S)$.
Similarly, for a vertex $i \in \bar S$, we denote its degree at
$G[\bar S]$ by $d_i(\bar S)$. We denote the maximum degrees, the
minimum degrees, the average degrees of $G[S]$ and $G[\bar S]$ by
$\Delta(S), \Delta(\bar S), \delta(S), \delta(\bar S), \bar d(S)$
and $\bar d(\bar S)$ respectively.   For an edge $e_p \in E(S, \bar
S)$, $t(e_p)$ of its vertices are in $S$, where $1 \le t(e_p) \le
k-1$. For all edges $e_p \in E(S, \bar S)$, the average value of
such $t(e_p)$ is denoted $t(S)$. Then $1 \le t(S) \le k-1$.
Similarly, we may define $t(\bar S)$. Then $t(S) + t(\bar S) = k$.
We call the minimum or maximum cardinality of such an edge cut the
{\bf edge connectivity} or {\bf maximum cut} of $G$, and denote it
by $e(G)$ or $c(G)$ respectively.

For $e_p = (i_1, \ldots, i_k) \in E$, define a $k$th order
$n$-dimensional symmetric tensor $\Q(e_p)$ by
$$\Q(e_p)x^k = \sum_{j=1}^k x_{i_j}^k + k x_{i_1}\cdots x_{i_k}$$
for any $x \in C^n$.   Then, for any $x \in C^n$, we have
$$\Q x^k = \sum_{e_p \in E} \Q(e_p)x^k.$$

\begin{Proposition} \label{p1}   For a $k$-graph $G$, we have the
following conclusions.

(a). The edge connectivity satisfies $e(G) \le \delta$.

(b). We have
$$c(G) \le {n \over k}(2\bar d - \delta).$$

(c). If $n \le 2k-1$, then $e(G) = \delta$.
\end{Proposition}

\noindent {\bf Proof.}  (a).  Assume that $d_j = \delta$.  Let $S =
\{ j \}$.  Then $|E(S, \bar S)| = d_j = d_{min}$.  This proves (a).

(b). Let $S$ be a nonempty proper subset of $V$.   Let $x = {1 \over
|S|^{1 \over k}}\sum_{i \in S} \ii$.  For $e_p \in E(S)$, we have
$$\Q(e_p)x^k = {2k \over |S|}.$$
For $e_p \in E(\bar S)$, we have
$$\Q(e_p)x^k = 0.$$
For $e_p \in E(S, \bar S)$, we have
$$\Q(e_p)x^k = {t(e_p) \over |S|}.$$
As
$$\Q x^k = \left( \sum_{e_p \in E(S)} + \sum_{e_p \in E(\bar S)} +
\sum_{e_p \in E(S, \bar S)}\right)\Q(e_p)x^k,$$ we have
\begin{equation} \label{2.10}
\Q x^k = {2k \over |S|}|E(S)| + {t(S) \over |S|}|E(S, \bar S)|.
\end{equation}
Similarly, letting $y = {1 \over |\bar S|^{1 \over k}}\sum_{i \in
\bar S} \ii$, we have
\begin{equation} \label{2.11}
\Q y^k = {2k \over |\bar S|}|E(\bar S)| + {t(\bar S) \over |\bar
S|}|E(S, \bar S)|.
\end{equation}
By (\ref{e7.1}) and (\ref{2.10}), we have
\begin{equation} \label{2.12}
|S| \delta \le 2k |E(S)| + t(S)|E(S, \bar S)|.
\end{equation}
By (\ref{e7.1}) and (\ref{2.11}), we have
\begin{equation} \label{2.13}
|\bar S| \delta \le 2k |E(\bar S)| + t(\bar S)|E(S, \bar S)|.
\end{equation}
Summing (\ref{2.12}) and (\ref{2.13}), we have
$$n \delta \le 2k \left(|E(S)| + |E(\bar S)|\right) + k |E(S, \bar
S)|,$$ i.e.,
$$n \delta \le 2k \left(m - |E(S, \bar S)|\right) + k |E(S, \bar
S)|,$$ which implies that
$$\delta \le {2km \over n} - {k \over n}|E(S, \bar S)|.$$
Noticing that $\bar d = {km \over n}$, we have $$|E(S, \bar S)| \le
{n \over k}( 2\bar d - d_{min}).$$   This proves (b).

(c). When $n \le 2k-1$, either $|S| < k$ or $|\bar S| < k$. Without
loss of generality, assume that $|S| < k$.  Then $E(S) = \emptyset$
and $|E(S)| = 0$.   From (\ref{2.12}), we have
$$|S| \delta \le t(S)|E(S, \bar S)|.$$
We always have $t(S) \le |S|$.   Thus, we have
$$\delta \le |E(S, \bar S)|.$$   Combining this with Conclusion
(a), we have Conclusion (c). \ep

\section{Analytic Connectivity}

We define the {\bf analytic connectivity} $\alpha(G)$ of the
$k$-graph $G$ by
$$\alpha(G) = \min_{j=1,\cdots, n} \min \{ \LL x^k : x \in \Re^n_+, \sum_{i=1}^n x_i^k = 1, x_j = 0 \}.$$

By Theorem \ref{t8}, $\LL x^k \ge 0$ for any $x \in \Re^n_+$. Thus,
$\alpha(G) \ge 0$.   We first prove the following proposition.

\begin{Proposition}  \label{p2}
The $k$-graph $G$ is connected if and only if the algebraic
connectivity $\alpha(G) > 0$.
\end{Proposition}

\noindent {\bf Proof.} Suppose that $G$ is not connected.  Let $G_1
= (V_1, E_1)$ be a component of $G$.  Then there is a $j \in V
\setminus V_1$.     Let $x = {1 \over |V_1|}\sum_{i \in V_1} \ii$.
Then $x$ is a feasible point of $\min \{ \LL x^k : x \in \Re^n_+,
\sum_{i=1}^n x_i^k = 1, x_j = 0 \}$, and we see that $\min \{ \LL
x^k : x \in \Re^n_+, \sum_{i=1}^n x_i^k = 1, x_j = 0 \} = 0$. This
implies that $\alpha(G) = 0$.

Suppose that $\alpha(G) = 0$.  There is a $j$ such that $\min \{ \LL
x^k : x \in \Re^n_+, \sum_{i=1}^n x_i^k = 1, x_j = 0 \} = 0$.
Suppose that $x^*$ is a minimizer of this minimization problem. Then
$x_j^* = 0$, $\LL (x^*)^k = 0$ and by optimization theory, there is
a Lagrange multiplier $\mu$ such that for $i = 1, \ldots n, i \not =
j$, either $x_i^* = 0$ and
\begin{equation} \label{e8.1}
d_i(x_i^*)^{k-1} - \sum \left\{ {1 \over (k-1)!}x^*_{i_2}\cdots
x^*_{i_k} : (i, i_2, \ldots, i_k ) \in E \right\} \ge \mu
(x_i^*)^{k-1},
\end{equation}
or $x_i^* > 0$ and
\begin{equation} \label{e8.2}
d_i(x_i^*)^{k-1} - \sum \left\{ {1 \over (k-1)!}x^*_{i_2}\cdots
x^*_{i_k} : (i, i_2, \ldots, i_k ) \in E \right\} = \mu
(x_i^*)^{k-1}.
\end{equation}
In (\ref{e8.1}) and (\ref{e8.2}), we always have $x^* \in \Re^n_+$,
$\sum_{i=1}^n (x_i^*)^k = 1$ and $x_j^* = 0$. Multiplying
(\ref{e8.1}) and (\ref{e8.2}) with $x_i^*$ and summing them
together, we have $\mu \sum_{i=1}^n (x_i^*)^k = \LL (x^*)^k = 0$,
i.e., $\mu = 0$.   Then for $i = 1, \ldots n, i \not = j$, either
$x_i^* = 0$ or
\begin{equation} \label{e8.3}
d_i(x_i^*)^{k-1} - \sum \left\{ {1 \over (k-1)!}x^*_{i_2}\cdots
x^*_{i_k} : (i, i_2, \ldots, i_k ) \in E \right\} = 0.
\end{equation}
 Let $x_r^* =
\max \{ x_i^* : i = 1, \ldots, n \}$.  Then by (\ref{e8.3}), we have
$$0 = d_r - \sum \left\{ {1 \over (k-1)!}{x^*_{i_2} \over x_r^*} \cdots
{x^*_{i_k} \over x^*_r} : (r, i_2, \ldots, i_k ) \in E \right\}.$$
Note that
$$d_r = \sum \left\{ {1 \over (k-1)!} : (r, i_2, \ldots, i_k ) \in E \right\}.$$
Thus, we have $x_i = x_r$ as long as $i$ and $r$ are in the same
edge.  From this, we see that $x_i = x_r$ as long as $i$ and $r$ are
in the same component of $G$.  Since $x_j^* = 0$, we see that $j$
and $r$ are in the different components of $G$, i.e., $G$ is not
connected. This proves the proposition.
 \ep

We now further explore an application of $\alpha(G)$.

\begin{Proposition} \label{p3}   For a $k$-graph $G$, we have
$$e(G) \ge {n \over k}\alpha(G).$$
\end{Proposition}

\noindent {\bf Proof.}  Let $S$ be a nonempty proper subset of $V$.
Then there is a $j \not \in S$ such that
\begin{equation} \label{22.1}
\min \{ \LL x^k : x \in \Re^n_+, \sum_{i=1}^n x_i^k = 1, x_j = 0 \}
\ge \alpha(G).
\end{equation}
 Let $x = {1 \over |S|^{1 \over k}}\sum_{i \in S} \ii$. Then $x$ is a feasible point of the minimization
 problem in (\ref{22.1}).  For
$e_p \in E(S)$ and $e_p \in E(\bar S)$, we have
$$\LL (e_p)x^k = 0,$$
where $\LL (e_p)$ is defined in Section 4.    For $e_p \in E(S, \bar
S)$, we have
$$\LL(e_p)x^k = {t(e_p) \over |S|}.$$
As
$$\LL x^k = \left( \sum_{e_p \in E(S)} + \sum_{e_p \in E(\bar S)} +
\sum_{e_p \in E(S, \bar S)}\right)\LL (e_p)x^k,$$ we have
\begin{equation} \label{22.10}
\LL x^k = {t(S) \over |S|}|E(S, \bar S)|.
\end{equation}
Similarly, letting $y = {1 \over |\bar S|^{1 \over k}}\sum_{i \in
\bar S} \ii$, we have
\begin{equation} \label{22.11}
\LL y^k = {t(\bar S) \over |\bar S|}|E(S, \bar S)|.
\end{equation}
By (\ref{22.1}) and (\ref{22.10}), we have
\begin{equation} \label{22.12}
|S| \alpha(G) \le t(S)|E(S, \bar S)|.
\end{equation}
By (\ref{22.1}) and (\ref{22.11}), we have
\begin{equation} \label{22.13}
|\bar S| \alpha(G) \le t(\bar S)|E(S, \bar S)|.
\end{equation}
Summing up (\ref{22.12}) and (\ref{22.13}), we have
$$n \alpha(G) \le k |E(S, \bar
S)|,$$ i.e.,
$${n \over k}\alpha(G) \le |E(S, \bar S)|.$$
This implies that
$$e(G) \ge {n \over k}\alpha(G).$$
\ep

We now give an upper bound for $\alpha(G)$.

\begin{Proposition} \label{p4}   For a $k$-graph $G$, we have
$$0 \le \alpha(G) \le \delta.$$
\end{Proposition}

\noindent {\bf Proof.}  We know $\alpha(G) \ge 0$. It suffices to
prove that $\alpha(G) \le \delta$.   Suppose that $d_r = \delta$ and
$j \not = r$.   Then $l^{(r)}$ is a feasible point of
$$\min \{ \LL x^k : x \in \Re^n_+, \sum_{i=1}^n x_i^k = 1, x_j = 0
\},$$ and $\LL (l^{(r)})^k = \delta$.   This implies that
$$\alpha(G) \le \min \{ \LL x^k : x \in \Re^n_+, \sum_{i=1}^n x_i^k = 1, x_j = 0
\} \le \delta.$$ \ep

By Proposition \ref{p2}, when $G$ is not connected, $\alpha(G) = 0$.
Let $n=k, m=1$, $E = \{ (1, 2, \ldots, k) \}$.   Then $\LL x^k =
\sum_{i=1}^k x_i^k -kx_1\cdots x_k$, and we see that $\alpha(G) = 1
= \delta$. Thus, both the lower bound $0$ and the upper bound
$\delta$ in Proposition \ref{p4} are attainable. However, it is
possible that $0 < \alpha(G) < \delta$. Let $k=3, n=4, m = 2$, $E =
\{ (1, 2, 3), (2, 3, 4) \}$.  Then $G$ is connected and $\alpha(G)
> 0$.   We have $\LL x^3 = x_1^3 + 2 x_2^3 + 2 x_3^3 + x_4^3 -
3x_1x_2x_3 - 3x_2x_3x_4$.   Consider
$$\min \{ \LL x^3 : x \in \Re^4_+, \sum_{i=1}^4 x_i^3 = 1, x_4 = 0 \}$$
$$= \min \{ x_1^3 + 2 x_2^3 + 2 x_3^3 - 3x_1x_2x_3 : x_1^3 + x_2^3 + x_3^3 = 1, x_1, x_2, x_3 \ge 0 \}.$$
Let $y = \left({1 \over 3}\right)^{1 \over 3}(1, 1, 1)$.   Then we
see that
$$\alpha(G) \le \min \{ \LL x^3 : x \in \Re^4_+, \sum_{i=1}^4 x_i^3 = 1, x_4 = 0
\} \le y_1^3 + 2 y_2^3 + 2 y_3^3 - 3y_1y_2y_3 = {2 \over 3} < 1 =
\delta.$$ Actually, the exact value of $\alpha(G)$ for this example
is $\alpha(G) = 1 - \beta^2$, where $\beta$ satisfies $\beta +
\beta^3 = 1$ and $0.5 < \beta < 1$.

{\bf Question 3}.  In general, how can we calculate $\alpha(G)$?

\section{Final Remarks}

In this paper, we propose a simple and natural definition for the
Laplacian and the signless Laplacian tensors of a uniform
hypergraph. We show that they have very nice spectral properties.
This sets the base for further exploring their applications in
spectral hypergraph theory. Several further questions are raised. We
expect that the research on these two Laplacian tensors will also
motivate the further development of spectral theory of tensors. Some
very recent papers \cite{HQ1, HQ2, HQS, HQX, Qi2, QXX} demonstrated
the impacts on these two aspects.

{\bf Acknowledgment.} The authors are very grateful to the two
referees for their valuable suggestions and comments, which have
considerably improved the presentation of the paper.

%%%%%%%%%%%%%%%%%%%%%%%%%%%%%%%%%%%%%%%%%%%%%%%%%%%%%%%%%%%%%%%%%%%%%%%%%%%%%%%%%%%%%%%

\end{document}